\documentclass{amsart}

\usepackage{amssymb, hyperref}

\usepackage{cite}

\usepackage[capitalise]{cleveref}

\crefname{thm}{Thm.}{}
\crefname{prop}{Prop.}{}
\crefname{lem}{Lem.}{}
\crefname{cor}{Cor.}{}

\usepackage[msc-links]{amsrefs}

\hypersetup{
  colorlinks   = true, 
  urlcolor     = blue, 
  linkcolor    = blue, 
  citecolor   = red 
}

\newtheorem{thm}{Theorem}

\newtheorem {lem}{Lemma}

\newtheorem {exa}{Example}

\newtheorem {rem}{Remark}
\newtheorem {cor}{Corollary}

\DeclareMathOperator{\lcm}{lcm}

\newcommand\Z{\mathbb Z}
\newcommand\Q{\mathbb Q}

\newcommand\w{\mathfrak w}

\DeclareMathOperator\wgcd{wgcd \, }
\DeclareMathOperator\awgcd{\overline{wgcd}\, }

\def\l{\lambda}

\def\x{\mathbf x}

\title{Computing efficiently the weighted greatest common divisor}

\author{Orgest Zaka}
\address{Department of Mathematics and Informatics, Agricultural University of Tirana \\
Research Institute of Science and Technology, Vlore, Albania}
 
\email{ozaka@ubt.edu.al, ozaka@risat.org}



\begin{document}

\maketitle

\begin{abstract}
In this paper we included some basic properties for weighted greatest common divisors, and discuss how to speed up computing the weighted greatest common divisor. By ordering the 'weights' we are able to significantly shorten the operations to computing $\wgcd$. In the absence of an efficient algorithm for computing $\wgcd$ by ordering the weights,  and using $\gcd$, we significantly reduce the numbers for which we want to compute $\wgcd$. As a final result in this paper we prove that: If $\x=(x_{0},\dots ,x_{n})\in \Z^{n+1}$, with weights $\w=(q_{0},\dots ,q_{n})$ and $q_{0}\leq \cdots \leq q_{n}$, then
$\wgcd_w(\x)=\wgcd_w(y_0,y_1,\dots,y_n)$, where $y_i=\gcd(x_i,\dots,x_n)$, and $y_0\leq y_1 \leq\dots \leq y_n$.  
This results, is not an algorithm to calculate $\wgcd$, but in some cases it can help by reducing the calculations, and speeding up the finding of $\wgcd$.
\end{abstract}

\section{Introduction}
The arithmetic of weighted projective spaces has seen a flurry of activity in the last few years; see \cite{MR627828, MR2403559, MR1798982, bombieri, MR2852925, m-sh, igor, f-sh, beshaj-guest, height-1, deng, hosgood, beshaj-polak} among others.  Most of this activity is due to the introduction of height on a weighted projective space as in \cite{sh-h} and \cite{b-g-sh}, see also \cite{MR1747280, MR879909, MR708341, silv-book, MR3045344, kaplansky}.  The concept of the weighted greatest common divisor was instrumental in the definition of the weighted height.    In this paper we investigate some of the arithmetic properties of the weighted greatest common divisors.

Let $\x = (x_0, \dots x_n ) \in \Z^{n+1}$ be a tuple of integers, not all equal to zero.  Their greatest common divisor, denoted by $\gcd (x_0, \dots, x_n)$, is defined as the largest integer $d$ such that $d | x_i$, for all $i=0, \dots , n$.  

The concept of the weighted greatest common divisor of a tuple for the ring of integers $\Z$ was defined in \cite{sh-h}.   Let $q_0$, \dots, $q_n$ be positive integers.  A set of weights is called the ordered tuple  $\w=(q_0, \dots, q_n)$.   Denote by  $r=\gcd (q_0, \dots , q_n)$ the greatest common divisor of $q_0, \dots , q_n$.  A \textit{weighted integer tuple} is a tuple $\x = (x_0, \dots, x_n ) \in \Z^{n+1}$ such that to each coordinate $x_i$ is assigned the weight $q_i$.  We multiply weighted tuples by scalars $\lambda \in \Q$ via 
\[ \lambda \star (x_0, \dots , x_n) = \left( \l^{q_0} x_0, \dots , \l^{q_n} x_n   \right) \]
For an ordered   tuple of integers  $\x=(x_0, \dots, x_n) \in \Z^{n+1}$, whose coordinates are not all zero, the \emph{weighted greatest common divisor with respect to the set of weights} $\w$ is the largest integer $d$ such that 
\[ d^{q_i} \, \mid \, x_i, \; \; \text{for all }    i=0, \dots, n.\]
Such integer $d$ does exist for any tuple $\x = (x_0, \dots, x_n ) \in \Z^{n+1}$; see \cite{b-g-sh}.  We will denote by  $\wgcd (x_0, \dots, x_n) = \wgcd (\x)$.

%
%
%

\section{Weighted greatest common divisors}
%

%
For a weighted integer  tuple $\x=(x_0, \dots, x_n)$ with weights $\w=(q_0, \dots, q_n)$ such that $g= \prod_{i=1}^t p_i^{s_i}$ is the factorization intro primes  of $g = \gcd (x_0, \dots, x_n)$, the weighted greatest common divisor $d = \wgcd (\x)$  is given by
\[
d = \prod_{i=1}^t p_i^{\alpha_i}
\]
where  $\alpha_i$ are the largest integers such that $d^{q_i}$ divides $x_i$ and $\alpha_i \leq \beta_i$. With $\beta_i=\min \left\{ {\bar q}_{j,i}, j =0,\ldots,n \right\}$ where  ${\bar q}_{j,i}$'s are the quotients in the integer Euclidean divisions, i.e.,  
\[ s_i = q_j{\bar q}_{j,i} +r_{j,i}, \quad 0 \leq r_{j,i} < q_j, \quad \forall i, i= 1,\ldots,t, \quad \forall j, j= 0,\ldots,n. \]
%
Let us now denote by $m= \lcm (q_0, \dots , q_n)$. 
%
For a weighted integer  tuple $\x=(x_0, \dots , x_n)$ with weights $\w=(q_0, \dots , q_n)$ such that $m= \lcm (q_0, \dots , q_n)$, the weighted greatest common divisor $d=\wgcd (\x)$  is given by
\[
d^m = \gcd \left(x_0^{\frac m {q_0}}, \dots , x_n^{\frac m {q_n}} \right)
\]
A tuple $\x=(x_0, \dots , x_n)$ with $\wgcd (\x)=1$ is called \textbf{normalized}. 

\section{Properties of the weighted greatest common divisors}
Let $\x=(x_{0},\dots ,x_{n})\in \Z^{n+1}$ be a tuple of integers, not all equal to zero, with weights $\w=(q_{0}, \ldots , q_{n})$.   Let the factorization of the integers $x_i$,  $(i=0, \ldots,  n)$ into
primes:
\[
x_{i}=\prod_{j=1}^{t}p_{j}^{\alpha _{j,i}},\text{ \ \ \ }\alpha _{j,i}\geq 0,j=1,\ldots,t
\]
Hence,  $\x$ is written as 
\[
\x=\left( \prod_{j=1}^{t}p_{j}^{\alpha _{j,0}}, \ldots , \prod_{j=1}^{t}p_{j}^{\alpha_{j,n}}\right)
\]
Then the greatest common divisor $g=\gcd (x)=\gcd (x_{0},\ldots x_{n})$ is given by
\[
g=\gcd (\x)=\gcd (x_{0},\ldots x_{n})=\prod_{j=1}^{t}p_{j}^{\min \{\alpha _{j,i}: i=0, \ldots, n\}}
\]
On the other hand, the weighted greatest common divisor $d=\wgcd(\x)$ is given by
\[
d= \wgcd (\x)=\prod_{j=1}^{t}p_{j}^{\min \left\{ \left\lfloor  \frac{\alpha _{j,i}}{q_{i}}\right\rfloor :i=0,\ldots,n\right\} }.
\]
It is obvious, that, if $\w=(1,\ldots ,1),$ well $q_{i}=1,\forall i=0,...,n,$ have 
\[
\wgcd(\x)=\gcd (\x),
\]
because $\ \min \{\alpha _{j,i}:i=0,\ldots ,n\}=\min \left\{ \left\lfloor  \frac{\alpha _{j,i}}{1}\right\rfloor :i=0,\ldots ,n\right\}$.

For a tuple $\x = ( x_0, \ldots , x_n ) \in \Z^{n+1}$ and a set of weights $\w=(q_{0}, \ldots, q_{n})$   we will denote by 
\[ 
\x^i := (x_0, \dots , x_{i-1}, x_{i+i}, \ldots , x_n) \in \Z^n,  \quad \text{ and } \quad \w^i=(q_{0}, \ldots, q_{i-1}, q_{i+1},\ldots, q_{n})
\]
for each $i=0, \ldots , n$.   Then we have the following. 
 
\begin{thm} \label{thm.1}
For $\x=(x_{0},\ldots,x_{n})\in \Z^{n+1},$ with weights $\w=(q_{0},\ldots,q_{n})$,  the following holds 
\begin{equation}
\wgcd_\w ( \x ) = \wgcd_{(1, q_i)}  ( \wgcd_{\w^i}  (\x^i) , x_i  ), 
\end{equation}
for each $i=0, \ldots , n$. 
\end{thm}

\proof

Rewrite $\x$ into primes:
\[
\x=\left( \prod_{j=1}^{t}p_{j}^{\alpha_{j,0}},\ldots,\prod_{j=1}^{t}p_{j}^{\alpha_{j,i}},\ldots,\prod_{j=1}^{t}p_{j}^{\alpha _{j,1}}\right) 
\]
we have that,
\[
\wgcd_\w ( \x )=\prod_{j=1}^{t}p_{j}^{\min \left\{ \left\lfloor \frac{\alpha _{j,k}}{q_{k}}\right\rfloor :k=0,\ldots,n\right\} }.
\]
On the other hand, have $\x^i := (x_, \ldots , x_{i-1}, x_{i+i}, \ldots , x_n) \in \Z^n$, with weights $\w^i=(q_{0}, \ldots, q_{i-1}, q_{i+1}, \ldots,q_{n})$. Rewrite now $\x^i$ into primes:
\[
\x^i=\left( \prod_{j=1}^{t}p_{j}^{\alpha_{j,0}},\ldots,\prod_{j=1}^{t}p_{j}^{\alpha_{j,i-1}},\prod_{j=1}^{t}p_{j}^{\alpha_{j,i+1}},\ldots,\prod_{j=1}^{t}p_{j}^{\alpha _{j,1}}\right) 
\]
for this, we have
\[
\wgcd_{\w^i}( \x^i )=\prod_{j=1}^{t}p_{j}^{\alpha_{j}^{\prime}}.
\]
where, $\alpha_{j}^{\prime}=\min \left\{ \left\lfloor  \frac{\alpha _{j,0}}{q_{0}}\right\rfloor,\ldots,\left\lfloor  \frac{\alpha _{j,i-1}}{q_{i-1}}\right\rfloor,\left\lfloor  \frac{\alpha _{j,i+1}}{q_{i+1}}\right\rfloor, \ldots, \left\lfloor 
\frac{\alpha _{j,n}}{q_{n}}\right\rfloor\right\}$.
Now write the ordered pairs $(\wgcd_{\w^i}( \x^i ),x_i)$ into primes,
\[
(\wgcd_{\w^i}( \x^i ),x_i)=\left(\prod_{j=1}^{t}p_{j}^{\alpha_{j}^{\prime}},\prod_{j=1}^{t}p_{j}^{\alpha_{j,i}}\right),
\]
and we have 
\[
\wgcd_{(1,q_i)} \left(\wgcd_{\w^i}( \x^i ),x_i\right) =\prod_{j=1}^{t}p_{j}^{\alpha_{j}^{\prime\prime}}.
\]
where  $ \alpha_{j}^{\prime\prime}=\min{\left\{\left\lfloor \frac{\alpha_{j}^{\prime}}{1}\right\rfloor,\left\lfloor \frac{\alpha_{j,i}}{q_{i}}\right\rfloor\right\}}$.
Therefore, 
\[
\begin{split}
\alpha_{j}^{\prime\prime}&=\min{\left\{\left\lfloor \frac{\alpha_{j}^{\prime}}{1}\right\rfloor,\left\lfloor \frac{\alpha_{j,i}}{q_{i}}\right\rfloor\right\}} \\
&=\min{\left\{\left\lfloor \frac{\min \left\{ \left\lfloor 
\frac{\alpha _{j,0}}{q_{0}}\right\rfloor,\ldots,\left\lfloor 
\frac{\alpha _{j,i-1}}{q_{i-1}}\right\rfloor,\left\lfloor 
\frac{\alpha _{j,i+1}}{q_{i+1}}\right\rfloor, \ldots, \left\lfloor 
\frac{\alpha _{j,n}}{q_{n}}\right\rfloor\right\}}{1}\right\rfloor,\left\lfloor \frac{\alpha_{j,i}}{q_{i}}\right\rfloor\right\}} \\
&=\min{\left\{\min \left\{ \left\lfloor 
\frac{\alpha _{j,0}}{q_{0}}\right\rfloor,\ldots,\left\lfloor 
\frac{\alpha _{j,i-1}}{q_{i-1}}\right\rfloor,\left\lfloor 
\frac{\alpha _{j,i+1}}{q_{i+1}}\right\rfloor, \ldots, \left\lfloor 
\frac{\alpha _{j,n}}{q_{n}}\right\rfloor\right\},\left\lfloor \frac{\alpha_{j,i}}{q_{i}}\right\rfloor\right\}} \\
&=\min{ 
\left\{\left\lfloor \frac{\alpha _{j,0}}{q_{0}}\right\rfloor,\ldots,
\left\lfloor \frac{\alpha _{j,i-1}}{q_{i-1}}\right\rfloor, 
\left\lfloor \frac{\alpha _{j,i+1}}{q_{i+1}}\right\rfloor, \ldots,  
\left\lfloor \frac{\alpha _{j,n}}{q_{n}}\right\rfloor, \left\lfloor \frac{\alpha_{j,i}}{q_{i}}\right\rfloor\right\}} \\
&=\min{ 
\left\{\left\lfloor \frac{\alpha _{j,0}}{q_{0}}\right\rfloor,\ldots,
\left\lfloor \frac{\alpha _{j,i-1}}{q_{i-1}}\right\rfloor, \left\lfloor \frac{\alpha_{j,i}}{q_{i}}\right\rfloor,
\left\lfloor \frac{\alpha _{j,i+1}}{q_{i+1}}\right\rfloor, \ldots,  
\left\lfloor \frac{\alpha _{j,n}}{q_{n}}\right\rfloor\right\}} \\
&=\min \left\{ \left\lfloor 
\frac{\alpha _{j,k}}{q_{k}}\right\rfloor :k=0,\ldots,n\right\} .
\end{split}
\]
Hence, 
\[
\begin{split}
\wgcd_{(1,q_i)} \left(\wgcd_{\w^i}( \x^i ),x_i\right) &=\prod_{j=1}^{t}p_{j}^{\alpha_{j}^{\prime\prime}} =\prod_{j=1}^{t}p_{j}^{\min \left\{ \left\lfloor  \frac{\alpha _{j,k}}{q_{k}}\right\rfloor :k=0,\ldots,n\right\}} = \wgcd_\w ( \x ).
\end{split}
\]
\qed

Let us  now see some examples.
\begin{exa}
Let $\x=(70352,5760,13824)$ and $\w=(2,2,3)$. Then 
\[
\wgcd_{\w}(70352,5760,13824) = 4.
\]
Following  \cref{thm.1} we have
\[
\wgcd_{\w}(70352,5760,13824)=\wgcd_{(2,1)}(70352,\wgcd_{(2,3)}(5760,13824))
\]
and 
\[
\wgcd_{(2,3)}(5760,13824) = \wgcd_{(2,3)}(2^7 \cdot 3^2 \cdot 5^1,2^9 \cdot 3^3 \cdot 5^0)= 24.
\]
Therefore, 
\[
\begin{split}
\wgcd_{\w}(70352,5760,13824)&= \wgcd_{(2,1)}(70352,\wgcd_{(2,3)}(5760,13824)) \\
&= \wgcd_{(2,1)}(70352,24) = 4.
\end{split}
\]
\end{exa}

Sometimes we can speed up the computation of the weighted greatest common divisor by the following result.

\begin{lem} \label{lem.1} 
For any integers $\x=(x_{0},x_{1})\in \Z^{2}$, with weights $%
\w=(q_{0},q_{1})$ \ we have

\begin{enumerate}
\item $\wgcd_{\w}(x_{0},x_{1})=\wgcd_{\w^{\prime }}(x_{1},x_{0}),$ (where $\w^{\prime }=(q_{1},q_{0})$),

\item $\wgcd_{\w}(x_{0},x_{1})=\wgcd_{\w}(\pm x_{0},\pm x_{1}),$ 

\item $\wgcd_{\w}(x_{0},x_{1})=\left\{ 
\begin{array}{c}
\wgcd_{\w}(x_{0},x_{1}-x_{0}),\text{ if \ }q_{0}\geq q_{1} \\ 
\wgcd_{\w}(x_{0}-x_{1},x_{1}),\text{ \ if  }q_{0}\leq q_{1}%
\end{array}%
\right. ,$ 
\end{enumerate}
\end{lem}

\proof
Points 1 and 2, are simple to verify. For the point 3, use factoriazations
into primes. For $(x_{0},x_{1})\in \Z^{2},$ have 
\[
(x_{0},x_{1})=\left( \prod_{j=1}^{t}p_{j}^{\alpha
_{j,0}},\prod_{j=1}^{t}p_{j}^{\alpha _{j,1}}\right) 
\]
and weighted greatest common divisor $\wgcd_{\w}(x_{0},x_{1})$ is given by%
\[
\wgcd_{\w}(x_{0},x_{1})=\prod_{j=1}^{t}p_{j}^{\min \left\{ \left\lfloor \frac{%
\alpha _{j,0}}{q_{0}}\right\rfloor ,\left\lfloor \frac{\alpha _{j,1}}{q_{1}}%
\right\rfloor \right\} }.
\]

\textbf{Case 1.} $q_{0}\geq q_{1}$. For $(x_{0},x_{1}-x_{0})\in \Z^{2},$ have 
\[
(x_{0},x_{1}-x_{0})=\left( \prod_{j=1}^{t}p_{j}^{\alpha
_{j,0}},\prod_{j=1}^{t}p_{j}^{\alpha _{j,1}}-\prod_{j=1}^{t}p_{j}^{\alpha
_{j,0}}\right) 
\]

Also we have,%
\[
\prod_{j=1}^{t}p_{j}^{\alpha _{j,1}}-\prod_{j=1}^{t}p_{j}^{\alpha
_{j,0}}=\prod_{j=1}^{t}p_{j}^{\min \{\alpha _{j,0},\alpha _{j,1}\}}\left[
\prod_{j=1}^{t}p_{j}^{\alpha _{j,1}-\min \{\alpha _{j,0},\alpha
_{j,1}\}}-\prod_{j=1}^{t}p_{j}^{\alpha _{j,0}-\min \{\alpha _{j,0},\alpha
_{j,1}\}}\right] 
\]

(in the factor $\left[ \prod_{j=1}^{t}p_{j}^{\alpha _{j,1}-\min \{\alpha
_{j,0},\alpha _{j,1}\}}-\prod_{j=1}^{t}p_{j}^{\alpha _{j,0}-\min \{\alpha
_{j,0},\alpha _{j,1}\}}\right] ,$ some of the primes are reduced, (are in 
zero-power)). Therefore,
\[
(x_{0},x_{1}-x_{0})=\left( \prod_{j=1}^{t}p_{j}^{\alpha
_{j,0}},\prod_{j=1}^{t}p_{j}^{\min \{\alpha _{j,0},\alpha _{j,1}\}}\left[
\prod_{j=1}^{t}p_{j}^{\alpha _{j,1}-\min \{\alpha _{j,0},\alpha
_{j,1}\}}-\prod_{j=1}^{t}p_{j}^{\alpha _{j,0}-\min \{\alpha _{j,0},\alpha
_{j,1}\}}\right] \right) 
\]

and we have%
\[
\begin{split}
\min \left\{ \left\lfloor \frac{\alpha _{j,0}}{q_{0}}\right\rfloor
,\left\lfloor \frac{\min \{\alpha _{j,0},\alpha _{j,1}\}}{q_{1}}%
\right\rfloor \right\}  &=\min \left\{ \left\lfloor \frac{\alpha _{j,0}}{%
q_{0}}\right\rfloor ,\min \left\{ \left\lfloor \frac{\alpha _{j,0}}{q_{1}}%
\right\rfloor ,\left\lfloor \frac{\alpha _{j,1}}{q_{1}}\right\rfloor
\right\} \right\}  \\
&=\min \left\{ \left\lfloor \frac{\alpha _{j,0}}{q_{0}}\right\rfloor
,\left\lfloor \frac{\alpha _{j,0}}{q_{1}}\right\rfloor ,\left\lfloor \frac{%
\alpha _{j,1}}{q_{1}}\right\rfloor \right\}  \\
&=\min \left\{ \left\lfloor \frac{\alpha _{j,0}}{q_{0}}\right\rfloor
,\left\lfloor \frac{\alpha _{j,1}}{q_{1}}\right\rfloor \right\} ,\left(
\left\lfloor \frac{\alpha _{j,0}}{q_{0}}\right\rfloor \leq \left\lfloor 
\frac{\alpha _{j,0}}{q_{1}}\right\rfloor \right) .
\end{split}
\]

Hence,%
\[
\wgcd_{\w}(x_{0},x_{1}-x_{0})=\prod_{j=1}^{t}p_{j}^{\min \left\{ \left\lfloor 
\frac{\alpha _{j,0}}{q_{0}}\right\rfloor ,\left\lfloor \frac{\alpha _{j,1}}{%
q_{1}}\right\rfloor \right\} }=\wgcd_{\w}\left( x_{0},x_{1}\right) 
\]

In the same way, we can write for

\textbf{Case 2. }$q_{0}\leq q_{1}$. For $(x_{0}-x_{1},x_{1})\in \Z^{2},$ have 
\[
(x_{0}-x_{1},x_{1})=\left( \prod_{j=1}^{t}p_{j}^{\alpha
_{j,0}}-\prod_{j=1}^{t}p_{j}^{\alpha _{j,1}},\prod_{j=1}^{t}p_{j}^{\alpha
_{j,1}}\right) 
\]

Also we have,%
\[
\prod_{j=1}^{t}p_{j}^{\alpha _{j,0}}-\prod_{j=1}^{t}p_{j}^{\alpha
_{j,1}}=\prod_{j=1}^{t}p_{j}^{\min \{\alpha _{j,0},\alpha _{j,1}\}}\left[
\prod_{j=1}^{t}p_{j}^{\alpha _{j,0}-\min \{\alpha _{j,0},\alpha
_{j,1}\}}-\prod_{j=1}^{t}p_{j}^{\alpha _{j,1}-\min \{\alpha _{j,0},\alpha
_{j,1}\}}\right] 
\]

(in the factor $\left[ \prod_{j=1}^{t}p_{j}^{\alpha _{j,0}-\min \{\alpha
_{j,0},\alpha _{j,1}\}}-\prod_{j=1}^{t}p_{j}^{\alpha _{j,1}-\min \{\alpha
_{j,0},\alpha _{j,1}\}}\right] ,$ some of the primes are reduced, (are in
zero-power)). Therefore,
\[
(x_{0}-x_{1},x_{1})=\left( \prod_{j=1}^{t}p_{j}^{\min \{\alpha _{j,0},\alpha
_{j,1}\}}\left[ \prod_{j=1}^{t}p_{j}^{\alpha _{j,0}-\min \{\alpha
_{j,0},\alpha _{j,1}\}}-\prod_{j=1}^{t}p_{j}^{\alpha _{j,1}-\min \{\alpha
_{j,0},\alpha _{j,1}\}}\right] ,\prod_{j=1}^{t}p_{j}^{\alpha _{j,1}}\right) 
\]
and we have%
\[
\begin{split}
\min \left\{ \left\lfloor \frac{\min \{\alpha _{j,0},\alpha _{j,1}\}}{q_{0}}%
\right\rfloor ,\left\lfloor \frac{\alpha _{j,1}}{q_{1}}\right\rfloor
\right\}  &= \min \left\{ \min \left\{ \left\lfloor \frac{\alpha _{j,0}}{%
q_{0}}\right\rfloor ,\left\lfloor \frac{\alpha _{j,1}}{q_{0}}\right\rfloor
\right\} ,\left\lfloor \frac{\alpha _{j,1}}{q_{1}}\right\rfloor \right\}  \\
&=\min \left\{ \left\lfloor \frac{\alpha _{j,0}}{q_{0}}\right\rfloor
,\left\lfloor \frac{\alpha _{j,1}}{q_{0}}\right\rfloor ,\left\lfloor \frac{%
\alpha _{j,1}}{q_{1}}\right\rfloor \right\}  \\
&=\min \left\{ \left\lfloor \frac{\alpha _{j,0}}{q_{0}}\right\rfloor
,\left\lfloor \frac{\alpha _{j,1}}{q_{1}}\right\rfloor \right\} ,\left(
\left\lfloor \frac{\alpha _{j,1}}{q_{1}}\right\rfloor \leq \left\lfloor 
\frac{\alpha _{j,1}}{q_{0}}\right\rfloor \right).
\end{split}
\]
Hence,%
\[
\wgcd_{\w}(x_{0}-x_{1},x_{1})=\prod_{j=1}^{t}p_{j}^{\min \left\{ \left\lfloor 
\frac{\alpha _{j,0}}{q_{0}}\right\rfloor ,\left\lfloor \frac{\alpha _{j,1}}{%
q_{1}}\right\rfloor \right\} }=\wgcd_{\w}\left( x_{0},x_{1}\right) 
\]
\qed

By point (2) of \cref{lem.1}, we have that: 
\[
\wgcd_{\w}(x_{0} -x_{1} ,x_{1})=\wgcd_{\w}(-[x_{0} -x_{1}] ,x_{1})=\wgcd_{\w}(x_{1} -x_{0} ,x_{1}).
\]
Hence,
\[
\wgcd_{\w}(x_{0},x_{1})=\wgcd_{\w}(x_{0} -x_{1} ,x_{1})=\wgcd_{\w}(x_{1} -x_{0} ,x_{1})
\]

\begin{cor}
For $(x_{0},x_{1})\in \Z^{2},$ with weights $\w=(q_{0},q_{1})$ that $q_{0}<q_{1}$,  have 
\[
\wgcd_{\w}(x_{0},x_{1})=\wgcd_{\w}\left( \gcd [x_{0},x_{1}],x_{1}\right) 
\]
\end{cor}

\proof
For $(x_{0},x_{1})\in \Z^{2},$ with weights $\w%
=(q_{0},q_{1}),$ and $q_{0}<q_{1}$ have 
\[
(x_{0},x_{1})=\left( \prod_{j=1}^{t}p_{j}^{\alpha
_{j,0}},\prod_{j=1}^{t}p_{j}^{\alpha _{j,1}}\right) 
\]
and weighted greatest common divisor $\wgcd_{\w}(x_{0},x_{1})$ is given by
\[
\wgcd_{\w}(x_{0},x_{1})=\prod_{j=1}^{t}p_{j}^{\min \left\{ \left\lfloor \frac{\alpha _{j,0}}{q_{0}}\right\rfloor ,\left\lfloor \frac{\alpha _{j,1}}{q_{1}} \right\rfloor \right\} }.
\]
On the other hand
\[
\gcd [x_{0},x_{1}]=\prod_{j=1}^{t}p_{j}^{\min \left\{ \alpha _{j,0},\alpha
_{j,1}\right\} }
\]
and
\[
\begin{split}
\wgcd_{\w}(\gcd [x_{0},x_{1}],x_{1}) &= \prod_{j=1}^{t}p_{j}^{\min \left\{ \left\lfloor \frac{\min \left\{ \alpha _{j,0},\alpha _{j,1}\right\} }{q_{0}}\right\rfloor ,\left\lfloor \frac{\alpha _{j,1}}{q_{1}}\right\rfloor \right\} } \\
&=\prod_{j=1}^{t}p_{j}^{\min \left\{ \left\lfloor \frac{\alpha _{j,0}}{q_{0}}\right\rfloor ,\left\lfloor \frac{\alpha _{j,1}}{q_{0}}\right\rfloor, \left\lfloor \frac{\alpha _{j,1}}{q_{1}}\right\rfloor \right\} } \\
&=\prod_{j=1}^{t}p_{j}^{\min \left\{ \left\lfloor \frac{\alpha _{j,0}}{q_{0}}\right\rfloor ,\left\lfloor \frac{\alpha _{j,1}}{q_{1}}\right\rfloor\right\} } \\
&= \wgcd_{\w}(x_{0},x_{1}).
\end{split}
\]
\qed

\begin{cor}   By point (2) of this \cref{lem.1} the following holds
\[
\wgcd_{\w}\left( x_{0},x_{1}\right) =\wgcd_{\w}\left( \left\vert x_{0}\right\vert, \left\vert x_{1}\right\vert \right) 
\]
\end{cor}

Next we see a few examples. 

\begin{exa}
Let $\x=(13824,5760)$ and $\w=(2,3)$. Then 
\[
\wgcd_{\w}(\x)=\wgcd_{\w}(13824-5760,5760)=\wgcd_{\w}(8064,5760)=4.
\]
Indeed,    $\x=(2^{9}\cdot 3^{3}\cdot 5^{0},2^{7}\cdot 3^{2}\cdot 5)$,  then
\[
\wgcd_{\w}(\x) = 2^{\min \left\{ \frac{9}{2},\frac{7}{3}\right\} }\cdot3^{^{\min \left\{ \frac{3}{2},\frac{2}{3}\right\} }}\cdot 5^{^{\min \left\{ \frac{0}{2},\frac{1}{3}\right\} }} =4.
\]
On the other hand,  $(8064,5760)=\left( 2^{7}\cdot 3^{2}\cdot 5^{0}\cdot
7,2^{7}\cdot 3^{2}\cdot 5\cdot 7^{0}\right)$, hence,
\[
\wgcd_{\w}(8064,5760) =\wgcd_{\w}\left( 2^{7}\cdot 3^{2}\cdot 5^{0}\cdot 7^{1},2^{7}\cdot
3^{2}\cdot 5\cdot 7^{0}\right) =4.
\]
We repeat this procedure once again, and we have it,
\[
\wgcd_{\w}(8064,5760) = \wgcd_{\w}(8064-5760,5760)  = \wgcd_{\w}(2304,5760)=4,
\]
because
\[
\wgcd_{\w}(2304,5760) = \wgcd_{(2,3)}\left( 2^{8}\cdot 3^{2}\cdot 5^{0},2^{7}\cdot 3^{2}\cdot 5^{1}\right) =4.
\]
\end{exa}


\begin{exa}
Let $\x=(5760,13824)$ and $\w=(2,3)$. Then 
\[
\wgcd_{\w}(\x)=\wgcd_{\w}(5760,13824)=24.
\]

If we use the fact that 
\[
\wgcd_{\w}(x_{0},x_{1})=\wgcd_{\w}\left( \gcd [x_{0},x_{1}],x_{1}\right),
\]

we have
\[
\begin{split}
\wgcd_{\w}(5760,13824) &=\wgcd_{\w}\left( \gcd
[5760,13824],13824\right) \\
&= \wgcd_{\w}(1152,13824),
\end{split}
\]

and
\[
\wgcd_{\w}(1152,13824) = \wgcd_{\w}(2^{7}\cdot 3^{2}\cdot 5^{0},2^{9}\cdot 3^{3}\cdot 5^{1}) =24.
\]
\end{exa}

From \cref{lem.1} and its corollary, we obtain the following result.

\begin{thm} \label{thm.2}
Let's have $\x=(x_{0}, x_{1})\in \Z^2$, with $\w=(q_0,q_1)$ and $q_{0} <q_{1}$.

\textbf{a)} If $x_{0} > x_{1}$, then 
\[ 
\wgcd (x_0, x_1) = \wgcd ( \alpha , x_1 ), 
\]
where $\alpha$ is the remainder when we divide $x_0$ by $x_1$.  In other words 
\[ 
x_0 = r \cdot x_1 + \alpha, 
\]
and $\alpha < x_1$.  So we have to deal now with smaller numbers $\alpha$ and $x_1$.  

\textbf{b)}  If $x_{1} > x_{0}$, then

\[ 
\wgcd (x_0, x_1) = \wgcd ( \beta , x_1 ), 
\]

where $\beta$ is the remainder when we divide $x_1$ by $x_0$.  In other words 
\[ 
x_1 = k \cdot x_0 + \beta, 
\]
and $\beta < x_0$, of course we have and $\beta < x_1$.  So we have to deal now with smaller numbers $\beta$ and $x_1$.  
\end{thm}

\proof 
\textbf{a)} We have that $x_0 > x_1$. From point (3) of \cref{lem.1}, we have:
\[  
\wgcd_{\w} (x_0, x_1) = \wgcd_{\w} (x_0 - x_1, x_1) 
\]
By using this repeatedly we can have 
\[
\begin{split}  
\wgcd_{\w} (x_0, x_1) &= \wgcd_{\w} (x_0 - x_1, x_1) \\
&= \wgcd_{\w} (x_0 - 2 \cdot x_1, x_1) \\
&= \wgcd_{\w} (x_0 - 3 \cdot x_1, x_1) \\
&  \vdots \\
&= \wgcd_{\w} (x_0 - r \cdot x_1, x_1)
\end{split}
\]
and we mark $ x_0 - r \cdot x_1 = \alpha $, where $\alpha$ is the remainder when we divide $x_0$ by $x_1$ and $\alpha < x_1$.  So we have to deal now with smaller numbers $\alpha$ and $x_1$. So
\[ 
\wgcd (x_0, x_1) = \wgcd ( \alpha , x_1 )
\]

\textbf{b)} We have that $x_0 < x_1$. From points (2) and (3) of \cref{lem.1}, we have:
\[
\begin{split}  
\wgcd_{\w} (x_0, x_1) &= \wgcd_{\w} (x_0 - x_1, x_1) \\
&=\wgcd_{\w}(-[x_{1} -x_{0}] ,x_{1}) \\
&=\wgcd_{\w}(x_{1} -x_{0},x_{1})
\end{split}
\]
By using this repeatedly we can have 
\[
\begin{split}  
\wgcd_{\w} (x_0, x_1) &=\wgcd_{\w}(x_{1} -x_{0},x_{1}) \\
&=\wgcd_{\w}(x_{1} -2 \cdot x_{0},x_{1}) \\
&=\wgcd_{\w}(x_{1} -3 \cdot x_{0},x_{1}) \\
& \vdots \\
&=\wgcd_{\w}(x_{1} - k \cdot x_{0},x_{1})
\end{split}
\]
and we mark $x_1 - k \cdot x_0 = \beta $, where $\beta$ is the remainder when we divide $x_1$ by $x_0$ and $\beta < x_0$.  So we have to deal now with smaller numbers $\beta$ and $x_1$.  So
\[ 
\wgcd (x_0, x_1) = \wgcd ( \beta , x_1 ), 
\]
\qed

\begin{exa}
Let $\x=(5760,13824)$ and $\w=(2,3)$. Then 
find the remainder of division, following case b) of \cref{lem.1}, ie
\[
13824=2\cdot 5760 +2304,
\]
from, \cref{lem.1}, we have that:
\[
\wgcd_{\w}(\x)=\wgcd(5760,13824)=\wgcd_{\w}(2304,13824)=24.
\]
Indeed, $\x=(2^{7}\cdot 3^{2}\cdot 5^{1},2^{9}\cdot 3^{3}\cdot 5^{0})$, then 
\[
\wgcd_{\w}(\x)=24.
\]
On the other hand,
\[
\wgcd_{\w}(2304,13824) = \wgcd_{\w}(2^{8}\cdot 3^{2}\cdot 5^{0},2^{9}\cdot 3^{3}\cdot 5^{0})=24.
\]
\end{exa}

\begin{exa}
Let $\x=(70352,13824)$ and $\w=(2,3)$. Then 
find the remainder of division, following case a) of \cref{lem.1}, ie
\[
70352=5\cdot 5760 +1232,
\]
from, \cref{lem.1}, we have that:
\[
\wgcd_{\w}(\x)=\wgcd_{\w}(70352,13824)=\wgcd_{\w}(1232,13824)=4.
\]
Calculation first
\[
\wgcd_{\w}(\x) = \wgcd_{\w}(2^{4}\cdot \cdot \cdot 4397^{1},2^{9} \cdot 3^{3} \cdot 5^{0}\cdots 4397^{0})=4.
\]
Mark with $\tilde{\x}=(1232,13824)$, now we calculate, $\wgcd_{\w}(\tilde{\x})$, and have

\[
\wgcd_{\w}(\tilde{\x})=\wgcd_{\w}(2^{4}\cdot 3^{0}\cdot 5^{0}\cdot 7^{1}\cdot 11^{1}, 2^{9}\cdot 3^{3}\cdot 5^{0}\cdot 7^{0}\cdot 11^{0})=4.
\]
\end{exa}

\begin{exa}
Let $\x=(70352,5760,13824)$ and $\w=(2,2,3)$. Then 
find the remainder of division, following case b) of \cref{lem.1}, ie
$$
13824=2\cdot 5760 +2304 \text{ and }\   70352=5\cdot 5760 +1232,
$$
now we have that,
\[
\wgcd_{\w}(70352,5760,13824)=\wgcd_{\w}(1232,2304,13824) = 4.
\]
\end{exa}

We can speed up the computation of the weighted greatest common divisor by the following result.

\begin{thm} \label{thm.3}
Let $\x=(x_{0},\dots ,x_{n})\in \Z^{n+1}$ be a tuple of
integers, not all equal to zero, with weights $\w=(q_{0},\ldots 
,q_{n})$, for which we have that $q_{0}<\dots <q_{n}$. Then we have:
\[
\wgcd_{\w}(\x)=\wgcd_{\w}\left( \gcd [x_{0},x_{1},\dots ,x_{n}],x_{1},\dots ,x_{n}\right) .
\]
\end{thm}

\proof Rewrite $\x $ into primes:
\[
\x= \left(\prod_{j=1}^{t}p_{j}^{\alpha _{j,0}},\prod_{j=1}^{t}p_{j}^{\alpha _{j,1}},\dots \prod_{j=1}^{t}p_{j}^{\alpha _{j,n}}\right).
\]%
Then we have
\[
\gcd(\x)=\prod_{j=1}^{t}p_{j}^{\min\left\{\alpha _{j,0},\alpha _{j,1},\dots ,\alpha _{j,n}\right\}}.
\]
Calculate now, $\wgcd\left( \gcd [x_{0},x_{1},\dots ,x_{n}],x_{1},\dots ,x_{n}\right)$,
\[
\begin{split}
\wgcd\left( \gcd [x_{0},x_{1},\dots ,x_{n}],x_{1},\dots ,x_{n}\right) &= \prod_{j=1}^{t}p_{j}^{\min\left\{\left\lfloor \frac{\min\left\{\alpha _{j,0},\alpha _{j,1},\dots ,\alpha _{j,n}\right\}}{q_0}\right\rfloor,\left\lfloor \frac{\alpha_{j,1}}{q_{1}}\right\rfloor ,\ldots , \left\lfloor \frac{\alpha_{j,n}}{q_n}\right\rfloor\right\}} \\
&= \prod_{j=1}^{t}p_{j}^{\min\left\{\left\lfloor \frac{\alpha_{j,0}}{q_0}\right\rfloor,\left\lfloor \frac{\alpha_{j,1}}{q_0}\right\rfloor,\ldots\left\lfloor \frac{\alpha_{j,n}}{q_0}\right\rfloor, \left\lfloor \frac{\alpha_{j,1}}{q_1}\right\rfloor,\ldots,\left\lfloor \frac{\alpha_{j,n}}{q_n}\right\rfloor\right\}} \\
&= \prod_{j=1}^{t}p_{j}^{\min\left\{\left\lfloor \frac{\alpha_{j,0}}{q_0}\right\rfloor,\left\lfloor \frac{\alpha_{j,1}}{q_1}\right\rfloor,\ldots,\left\lfloor \frac{\alpha_{j,n}}{q_n}\right\rfloor\right\}} \\
&= \wgcd\left(x_{0},x_{1},\dots ,x_{n}\right),
\end{split}
\]
in the above equations we use the fact that $\frac{\alpha_{j,i}}{q_0}<\frac{\alpha_{j,i}}{q_{i}}$, for $i=1,2,\ldots,n$.

Hence, 
\[
\wgcd\left( \gcd [x_{0},x_{1},\dots ,x_{n}],x_{1},\dots ,x_{n}\right)=\wgcd\left(x_{0},x_{1},\dots ,x_{n}\right).
\]

\qed 

As a generalization, for \cref{thm.3}, and a way to speed up the computation of the weighted greatest common divisor, we have the following result. 

\begin{thm} \label{myTh}
Let $\x \in \Z^{n+1}$, such that $\x= (x_0, \ldots , x_n)$ and    $\w=(q_{0}, \ldots  , q_{n})$, for which we have that $q_{0}<\dots <q_{n}$. 
Then  
\[
\wgcd_{\w}  ( \x) =\wgcd_{\w} \left(\gcd [x_{0},\dots ,x_{n}], \gcd[x_{1},\dots ,x_{n}],\ldots,\gcd[x_{n-1},x_{n}], x_{n} \right).
\]
\end{thm} 

\proof We act in the same way as the proof of \cref{thm.3}. Rewrite $\x$ into primes:
\[
\x= \left(\prod_{j=1}^{t}p_{j}^{\alpha _{j,0}},\prod_{j=1}^{t}p_{j}^{\alpha _{j,1}},\dots \prod_{j=1}^{t}p_{j}^{\alpha _{j,n}}\right).
\]%

Mark with,
\[
\begin{split}
y_0 &=\gcd\left(x_{0},\dots ,x_{n}\right)=\prod_{j=1}^{t}p_{j}^{\min\left\{\alpha _{j,0},\dots ,\alpha _{j,n}\right\}},\\
y_1 &=\gcd\left(x_{1},\dots ,x_{n}\right)=\prod_{j=1}^{t}p_{j}^{\min\left\{\alpha _{j,1},\dots ,\alpha _{j,n}\right\}},\\
y_2 &=\gcd\left(x_{2},\dots ,x_{n}\right)=\prod_{j=1}^{t}p_{j}^{\min\left\{\alpha _{j,2},\dots ,\alpha _{j,n}\right\}},\\
&\vdots \\
y_{n-1} &=\gcd\left(x_{n-1},x_{n}\right)=\prod_{j=1}^{t}p_{j}^{\min\left\{\alpha _{j,n-1},\alpha _{j,n}\right\}} \\
y_n &=x_n.
\end{split}
\]
Calculate now, $\wgcd_{\w}(y_0, y_1,\ldots, y_{n})$,
\[
\begin{split}
d&=\wgcd_{\w}(y_0, y_1,\ldots, x_{n}) \\
&= \prod_{j=1}^{t}p_{j}^{\min\left\{\left\lfloor \frac{{\min\left\{\alpha _{j,0},\alpha _{j,1},\dots ,\alpha _{j,n}\right\}}}{q_0}\right\rfloor; \left\lfloor \frac{{\min\left\{\alpha _{j,1},\alpha _{j,2},\dots ,\alpha _{j,n}\right\}}}{q_1}\right\rfloor;
 \ldots; \left\lfloor 
\frac{{\min\left\{\alpha _{j,n-1},\alpha _{j,n}\right\}}}{q_n-1}\right\rfloor;\left\lfloor \frac{\alpha_{j,n}}{q_n}\right\rfloor\right\}} \\
&= \prod_{j=1}^{t}p_{j}^{\min{\left\{\left\lfloor \frac{\alpha _{j,0}}{q_0}\right\rfloor,\ldots, \left\lfloor \frac{\alpha _{j,n}}{q_0}\right\rfloor; \left\lfloor \frac{\alpha _{j,1}}{q_1}\right\rfloor,\ldots,\left\lfloor \frac{\alpha _{j,n}}{q_1}\right\rfloor;\ldots;\left\lfloor \frac{\alpha _{j,n-1}}{q_{n-1}}\right\rfloor,\left\lfloor \frac{\alpha _{j,n}}{q_{n-1}}\right\rfloor; \left\lfloor \frac{\alpha _{j,n}}{q_n}\right\rfloor\right\}}}.
\end{split}
\]
But, for a sorting order of weights, so, $q_{0}<q_{1}<q_{2}<\dots <q_{n-1}<q_{n}$, we have the following  inequalities:
\[
\begin{split}
\frac{\alpha_{j,i}}{q_0} &< \frac{\alpha_{j,i}}{q_{i}}, \text{ for } i=1,2,\ldots,n, \\
\frac{\alpha_{j,i}}{q_1} &< \frac{\alpha_{j,i}}{q_{i}}, \text{ for } i=2,3,\ldots,n,\\
\frac{\alpha_{j,i}}{q_2} &< \frac{\alpha_{j,i}}{q_{i}}, \text{ for } i=3,4,\ldots,n,\\
&\vdots \\
\frac{\alpha_{j,i}}{q_{n-1}} &< \frac{\alpha_{j,i}}{q_{i}}, \text{ for } i=n.
\end{split}
\]
exploiting these inequalities we have that, 

\[
d = \prod_{j=1}^{t}p_{j}^{\min{\left\{\left\lfloor \frac{\alpha _{j,0}}{q_0}\right\rfloor; \left\lfloor \frac{\alpha _{j,1}}{q_1}\right\rfloor;\ldots;\left\lfloor \frac{\alpha _{j,n-1}}{q_{n-1}}\right\rfloor; \left\lfloor \frac{\alpha _{j,n}}{q_n}\right\rfloor\right\}}} =  \wgcd \left(x_{0},\dots ,x_{n}\right),
\]
thus,
\[ 
\wgcd_{\w}(y_0, y_1,\ldots, y_{n})=\wgcd_{\w}(x_0, x_1,\ldots, x_{n})
\]
Hence,
\[
\wgcd \left(x_{0},\dots ,x_{n}\right)=\wgcd\left(\gcd [x_{0},\dots ,x_{n}], \gcd[x_{1},\dots ,x_{n}],\ldots,\gcd[x_{n-1},x_{n}], x_{n} \right) .
\]
\qed

Given the fact that
\[
\gcd(a_0,a_1,...,a_n)\leq\gcd(a_1,...,a_n) \leq \gcd(a_2,...,a_n)\leq \ldots \leq {a_n}, 
\]
we have the following:

\begin{cor}
Let $\x=(x_{0},\dots ,x_{n})\in \Z^{n+1}$,   $\w=(q_{0},\ldots  , q_{n})$, such  that $q_{0}<\dots <q_{n}$. Then  
\[
\wgcd_{\w}(\x) =\wgcd_\w  \left(y_{0},\dots ,y_{n}\right)
\]
where $y_i = \gcd (x_i, \ldots , x_n)$ for $i=0, 1, \ldots , n$.  Moreover,    
\[
y_{0}\leq y_{1} \leq \cdots \leq y_{n}
\]
\end{cor}

This Theorem \ref{myTh}, it is not an algorithm to calculate weighted greatest common divisors, but in some cases it can help by reducing the calculations, and speeding up the finding of $\wgcd$. 

Below we give two examples, the first example when this theorem is not efficient, and another example when this theorem is efficient.

\begin{rem}
This Theorem \ref{myTh}, it is weak in the case, when we have the numbers, we have the power of prime numbers, but we really hope it helps in other cases, different from these.
\end{rem}
\begin{exa}
Consider weights $\w=(2,3)$ and a tuple $x=(p^2, p^3)$ for a prime $p$.  The weighted gcd of $x$ is 
\[
\wgcd (x) = p.
\]

According to the main theorem \ref{myTh} we have,
\[
\wgcd (p^2, p^3) = \wgcd (\gcd(p^2, p^3), p^3) = \wgcd (p^2, p^3),
\]
nothing really happens here.
\end{exa}

Two examples where the 'power' of my theorem is felt

\begin{exa}
Let $\x=(70352,5760,13824)$ and $\w=(2,2,3)$. Then 
\[
\wgcd_{\w}(70352,5760,13824) = 4.
\]
Following  \cref{myTh} we have
\[
\begin{split}
\wgcd_{\w}(70352,5760,13824)&=\wgcd_{\w}(\gcd[70352,5760,13824],\gcd[5760,13824],13824)\\
&=\wgcd_{\w}(16,1152,13824)\\
&=\wgcd_{\w}(2^4,2^7 \cdot 3^2,2^9 \cdot 3^3)=2^2=4.\\
\end{split}
\]
\end{exa}

\begin{exa}
Let us have, a tuple $\x=(123456, 243226, 5789534, 234566, 4322166)$, with weights $\w=(7,5,3,2,9)$. \\
We follow the main theorem \ref{myTh}, and we have that,
$\w'=(2,3,5,7,9)$, and 
\[ \x'=(234566, 5789534, 243226, 123456, 4322166),\]
so from Lemma \ref{lem.1} we have,
\[
\wgcd_{\w}(\x)=\wgcd_{\w'}(\x')
\]
From theorem \ref{myTh} have that,
\[
\begin{split}
\wgcd_{\w'}(\x') &=\wgcd_{\w'}(234566, 5789534, 243226, 123456, 4322166) \\
&=\wgcd_{\w'}(\gcd[234566, 5789534, 243226, 123456, 4322166], \\
&  \gcd[5789534, 243226, 123456, 4322166], \\
& \gcd[243226, 123456, 4322166], \gcd[123456, 4322166],  4322166)\\
&=\wgcd_{\w'}(2, 2, 2,6,4322166)=1.
\end{split}
\]
\end{exa}

\section{Declarations}

\subsection*{Conflict of Interest Statement}
There is no conflict of interest with any funder.



\bibliography{refs}

\begin{bibdiv}
\begin{biblist}

\bib{MR879909}{article}{
      author={Beltrametti, Mauro},
      author={Robbiano, Lorenzo},
       title={Introduction to the theory of weighted projective spaces},
        date={1986},
        ISSN={0723-0869},
     journal={Exposition. Math.},
      volume={4},
      number={2},
       pages={111\ndash 162},
      review={\MR{879909}},
}

\bib{beshaj-polak}{incollection}{
      author={Beshaj, L.},
      author={Polak, M.},
       title={On hyperelliptic curves of genus 3},
        date={2019},
   booktitle={Algebraic curves and their applications},
      series={Contemp. Math.},
      volume={724},
   publisher={Amer. Math. Soc., Providence, RI},
       pages={161\ndash 173},
         url={https://doi.org/10.1090/conm/724/14589},
      review={\MR{3916739}},
}

\bib{beshaj-guest}{incollection}{
      author={Beshaj, Lubjana},
      author={Guest, Scott},
       title={The weighted moduli space of binary sextics},
        date={2019},
   booktitle={Algebraic curves and their applications},
      series={Contemp. Math.},
      volume={724},
   publisher={Amer. Math. Soc., Providence, RI},
       pages={33\ndash 44},
         url={https://doi.org/10.1090/conm/724/14583},
      review={\MR{3916733}},
}

\bib{b-g-sh}{article}{
      author={Beshaj, Lubjana},
      author={Gutierrez, Jaime},
      author={Shaska, Tony},
       title={Weighted greatest common divisors and weighted heights},
        date={2019},
     journal={arXiv preprint arXiv:1902.06563},
}

\bib{MR2852925}{article}{
      author={Bini, Gilberto},
       title={Quotients of hypersurfaces in weighted projective space},
        date={2011},
        ISSN={1615-715X},
     journal={Adv. Geom.},
      volume={11},
      number={4},
       pages={653\ndash 667},
         url={https://doi.org/10.1515/ADVGEOM.2011.029},
      review={\MR{2852925}},
}

\bib{bombieri}{book}{
      author={Bombieri, Enrico},
      author={Gubler, Walter},
       title={Heights in {D}iophantine geometry},
      series={New Mathematical Monographs},
   publisher={Cambridge University Press, Cambridge},
        date={2006},
      volume={4},
        ISBN={978-0-521-84615-8; 0-521-84615-3},
         url={http://dx.doi.org/10.1017/CBO9780511542879},
      review={\MR{2216774 (2007a:11092)}},
}

\bib{MR627828}{article}{
      author={Buium, Alexandru},
       title={Weighted projective spaces as ample divisors},
        date={1981},
        ISSN={0035-3965},
     journal={Rev. Roumaine Math. Pures Appl.},
      volume={26},
      number={6},
       pages={833\ndash 842},
      review={\MR{627828}},
}

\bib{deng}{article}{
      author={Deng, An-Wen},
       title={Rational points of weighted projective spaces},
        date={1998},
     journal={arxiv:},
}

\bib{igor}{incollection}{
      author={Dolgachev, Igor},
       title={Weighted projective varieties},
        date={1982},
   booktitle={Group actions and vector fields ({V}ancouver, {B}.{C}., 1981)},
      series={Lecture Notes in Math.},
      volume={956},
   publisher={Springer, Berlin},
       pages={34\ndash 71},
         url={https://doi.org/10.1007/BFb0101508},
      review={\MR{704986}},
}

\bib{MR1747280}{incollection}{
      author={Fauntleroy, A.},
       title={Moduli of complete intersections in weighted projective spaces},
        date={1999},
   booktitle={African {A}mericans in mathematics, {II} ({H}ouston, {TX},
  1998)},
      series={Contemp. Math.},
      volume={252},
   publisher={Amer. Math. Soc., Providence, RI},
       pages={77\ndash 84},
         url={https://doi.org/10.1090/conm/252/1747280},
      review={\MR{1747280}},
}

\bib{f-sh}{incollection}{
      author={Frey, Gerhard},
      author={Shaska, Tony},
       title={Curves, {J}acobians, and cryptography},
        date={2019},
   booktitle={Algebraic curves and their applications},
      series={Contemp. Math.},
      volume={724},
   publisher={Amer. Math. Soc., Providence, RI},
       pages={279\ndash 344},
         url={https://doi.org/10.1090/conm/724/14596},
      review={\MR{3916746}},
}

\bib{MR3045344}{article}{
      author={Guo, Shuai},
      author={Zhou, Jian},
       title={Orbifold elliptic genera of quotients of complete intersections
  in weighted projective spaces},
        date={2013},
        ISSN={0129-167X},
     journal={Internat. J. Math.},
      volume={24},
      number={2},
       pages={1350012, 24},
         url={https://doi.org/10.1142/S0129167X13500122},
      review={\MR{3045344}},
}

\bib{silv-book}{book}{
      author={Hindry, Marc},
      author={Silverman, Joseph~H.},
       title={Diophantine geometry},
      series={Graduate Texts in Mathematics},
   publisher={Springer-Verlag, New York},
        date={2000},
      volume={201},
        ISBN={0-387-98975-7; 0-387-98981-1},
         url={http://dx.doi.org/10.1007/978-1-4612-1210-2},
        note={An introduction},
      review={\MR{1745599 (2001e:11058)}},
}

\bib{hosgood}{unpublished}{
      author={Hosgood, Timothy},
       title={An introduction to varieties in weighted projective space},
}

\bib{MR1798982}{incollection}{
      author={Iano-Fletcher, A.~R.},
       title={Working with weighted complete intersections},
        date={2000},
   booktitle={Explicit birational geometry of 3-folds},
      series={London Math. Soc. Lecture Note Ser.},
      volume={281},
   publisher={Cambridge Univ. Press, Cambridge},
       pages={101\ndash 173},
      review={\MR{1798982}},
}

\bib{kaplansky}{book}{
      author={Kaplansky, Irving},
       title={Commutative rings},
   publisher={The University of Chicago Press},
        date={1974},
}

\bib{m-sh}{article}{
      author={Malmendier, A.},
      author={Shaska, T.},
       title={From hyperelliptic to superelliptic curves},
        date={2019},
     journal={Albanian J. Math.},
      volume={13},
      number={1},
       pages={107\ndash 200},
      review={\MR{3978315}},
}

\bib{sh-h}{incollection}{
      author={Mandili, Jorgo},
      author={Shaska, Tony},
       title={Computing heights on weighted projective spaces},
        date={2019},
   booktitle={Algebraic curves and their applications},
      series={Contemp. Math.},
      volume={724},
   publisher={Amer. Math. Soc., Providence, RI},
       pages={149\ndash 160},
         url={https://doi.org/10.1090/conm/724/14588},
      review={\MR{3916738}},
}

\bib{MR2403559}{article}{
      author={Rossi, Carlo~A.},
       title={Weighted projective spaces and minimal nilpotent orbits},
        date={2008},
        ISSN={1088-4165},
     journal={Represent. Theory},
      volume={12},
       pages={208\ndash 224},
         url={https://doi.org/10.1090/S1088-4165-08-00328-2},
      review={\MR{2403559}},
}

\bib{height-1}{incollection}{
      author={Shaska, T.},
      author={Beshaj, L.},
       title={Heights on algebraic curves},
        date={2015},
   booktitle={Advances on superelliptic curves and their applications},
      series={NATO Sci. Peace Secur. Ser. D Inf. Commun. Secur.},
      volume={41},
   publisher={IOS, Amsterdam},
       pages={137\ndash 175},
      review={\MR{3525576}},
}

\bib{MR708341}{incollection}{
      author={Steenbrink, Joseph},
       title={On the {P}icard group of certain smooth surfaces in weighted
  projective spaces},
        date={1982},
   booktitle={Algebraic geometry ({L}a {R}\'abida, 1981)},
      series={Lecture Notes in Math.},
      volume={961},
   publisher={Springer, Berlin},
       pages={302\ndash 313},
         url={https://doi.org/10.1007/BFb0071290},
      review={\MR{708341}},
}

\end{biblist}
\end{bibdiv}


\begin{thebibliography}{99} 

\bibitem{beshaj-polak}Beshaj, L. \& Polak, M. On hyperelliptic curves of genus 3. {\em Algebraic Curves And Their Applications}. \textbf{724} pp. 161-173 (2019), https://doi.org/10.1090/conm/724/14589
\bibitem{beshaj-guest}Beshaj, L. \& Guest, S. The weighted moduli space of binary sextics. {\em Algebraic Curves And Their Applications}. \textbf{724} pp. 33-44 (2019), https://doi.org/10.1090/conm/724/14583
\bibitem{m-sh}Malmendier, A. \& Shaska, T. From hyperelliptic to superelliptic curves. {\em Albanian J. Math.}. \textbf{13}, 107-200 (2019)
\bibitem{f-sh}Frey, G. \& Shaska, T. Curves, Jacobians, and cryptography. {\em Algebraic Curves And Their Applications}. \textbf{724} pp. 279-344 (2019), https://doi.org/10.1090/conm/724/14596
\bibitem{sh-h}Mandili, J. \& Shaska, T. Computing heights on weighted projective spaces. {\em Algebraic Curves And Their Applications}. \textbf{724} pp. 149-160 (2019), https://doi.org/10.1090/conm/724/14588
\bibitem{b-g-sh}Beshaj, L., Gutierrez, J. \& Shaska, T. Weighted greatest common divisors and weighted heights. {\em ArXiv Preprint ArXiv:1902.06563}. (2019)
\bibitem{hosgood}Hosgood, T. An introduction to varieties in weighted projective space. 
\bibitem{kaplansky}Kaplansky, I. Commutative rings. (The University of Chicago Press,1974)
\bibitem{deng}Deng, A. Rational points of weighted projective spaces. {\em Arxiv:}. (1998)
\bibitem{height-1}Shaska, T. \& Beshaj, L. Heights on algebraic curves. {\em Advances On Superelliptic Curves And Their Applications}. \textbf{41} pp. 137-175 (2015)
\bibitem{MR1798982}Iano-Fletcher, A. Working with weighted complete intersections. {\em Explicit Birational Geometry Of 3-folds}. \textbf{281} pp. 101-173 (2000)
\bibitem{igor}Dolgachev, I. Weighted projective varieties. {\em Group Actions And Vector Fields (Vancouver, B.C., 1981)}. \textbf{956} pp. 34-71 (1982), https://doi.org/10.1007/BFb0101508
\bibitem{MR3045344}Guo, S. \& Zhou, J. Orbifold elliptic genera of quotients of complete intersections in weighted projective spaces. {\em Internat. J. Math.}. \textbf{24}, 1350012, 24 (2013), https://doi.org/10.1142/S0129167X13500122
\bibitem{MR2852925}Bini, G. Quotients of hypersurfaces in weighted projective space. {\em Adv. Geom.}. \textbf{11}, 653-667 (2011), https://doi.org/10.1515/ADVGEOM.2011.029
\bibitem{MR2403559}Rossi, C. Weighted projective spaces and minimal nilpotent orbits. {\em Represent. Theory}. \textbf{12} pp. 208-224 (2008), https://doi.org/10.1090/S1088-4165-08-00328-2
\bibitem{MR1747280}Fauntleroy, A. Moduli of complete intersections in weighted projective spaces. {\em African Americans In Mathematics, II (Houston, TX, 1998)}. \textbf{252} pp. 77-84 (1999), https://doi.org/10.1090/conm/252/1747280
\bibitem{MR879909}Beltrametti, M. \& Robbiano, L. Introduction to the theory of weighted projective spaces. {\em Exposition. Math.}. \textbf{4}, 111-162 (1986)
\bibitem{MR708341}Steenbrink, J. On the Picard group of certain smooth surfaces in weighted projective spaces. {\em Algebraic Geometry (La Rábida, 1981)}. \textbf{961} pp. 302-313 (1982), https://doi.org/10.1007/BFb0071290
\bibitem{MR627828}Buium, A. Weighted projective spaces as ample divisors. {\em Rev. Roumaine Math. Pures Appl.}. \textbf{26}, 833-842 (1981)
\bibitem{lang-2}Lang, S. Fundamentals of Diophantine geometry. (Springer-Verlag, New York,1983), http://dx.doi.org/10.1007/978-1-4757-1810-2
\bibitem{Fa}Faltings, G., Wüstholz, G., Grunewald, F., Schappacher, N. \& Stuhler, U. Rational points. (Friedr. Vieweg & Sohn, Braunschweig,1992), http://dx.doi.org/10.1007/978-3-322-80340-5, Papers from the seminar held at the Max-Planck-Institut für Mathematik, Bonn/Wuppertal, 1983/1984, With an appendix by Wüstholz
\bibitem{bombieri}Bombieri, E. \& Gubler, W. Heights in Diophantine geometry. (Cambridge University Press, Cambridge,2006), http://dx.doi.org/10.1017/CBO9780511542879
\bibitem{silv-book}Hindry, M. \& Silverman, J. Diophantine geometry. (Springer-Verlag, New York,2000), http://dx.doi.org/10.1007/978-1-4612-1210-2, An introduction
\end{thebibliography}

\end{document}